\documentclass{commat}

\title[Generalized curvature tensor and hypersurfaces of Kenmotsu Hermitian manifolds]{Generalized curvature tensor and the hypersurfaces of the Hermitian manifold for the class of Kenmotsu type}

\author{%
    Mohammed Y. Abass and Habeeb M. Abood
    }

\affiliation{
    \address{Mohammed Yousif Abass  --
    Department of Mathematics, College of Science, University of Basrah, Basrah, Iraq.
        }
    \email{%
    mohammedyousif42@yahoo.com
    }
    \address{Habeeb Mtashar Abood --
    Department of Mathematics, College of Education for Pure Sciences, University of Basrah, Basrah, Iraq
        }
    \email{%
    iraqsafwan2006@gmail.com
    }
    }

\abstract{%
    This paper determines the components of the generalized curvature tensor for the class of Kenmotsu type and establishes the mentioned class is $\eta$-Einstein manifold when the generalized curvature tensor is flat; the converse holds true under suitable conditions. It also introduces the notion of generalized $\Phi$-holomorphic sectional ($G\Phi SH$-) curvature tensor and thus finds the necessary and sufficient conditions for the class of Kenmotsu type to be of constant $G\Phi SH$-curvature. In addition, the notion of $\Phi$-generalized semi-symmetric is introduced and its relationship with the class of Kenmotsu type and $\eta$-Einstein manifold established. Furthermore, this paper generalizes the notion of the manifold of constant curvature and deduces its relationship with the aforementioned ideas. It finally shows that the class of Kenmotsu type exists as a hypersurface of the Hermitian manifold and derives a relation between the components of the Riemannian curvature tensors of the almost Hermitian manifold and its hypersurfaces.
    }



\VOLUME{32}
\NUMBER{1}
\YEAR{2024}
\firstpage{103}
\DOI{https://doi.org/10.46298/cm.10869}

\begin{paper}

\section{Introduction}\label{sec1}

The notion of generalized curvature tensor was introduced by Shaikh and Kundu \cite{SHK14} to generalize well-known curvature tensors such as the conformal curvature tensor, the concircular tensor, and the conharmonic tensor. Yildiz and De \cite{YID12} introduced and studied $\Phi$-projectively semisymmetric and $\Phi$-Weyl semisymmetric non-Sasakian $(k, \mu)$-contact metric manifolds while Kenmotsu \cite{Ken72} and Kirichenko and Khari-tonova \cite{KIK12} discussed the $\Phi$-holomorphic sectional curvature tensor. On the other hand, investigation of the geometry of the submanifolds of some Riemannian manifolds has captured the interest of authors such as Alegre and Carriazo \cite{ALC09}, Sular and \"{O}zg\"{u}r \cite{SUO09} and Chen \cite{Che17}.
The special subject in the study of the geometry of submanifolds is the hypersurface of the Riemannian manifolds, which has been discussed by Goldberg \cite{Gol68}. We concentrated on the geometry of the hypersurfaces of the almost Hermitian manifolds that have almost contact structures on the associated G-structure space. The last mentioned topic was studied by Banaru and Kirichenko \cite{BAK15}. Moreover, Ignatochkina \cite{Ign17}, Ignatochkina and Morozov \cite{IGM11}, and Nikiforova and Ignatochkina \cite{NII11} studied the transformations and conformal transformations on hypersurfaces induced from almost Hermitian manifolds.

\indent The aim of this article is organized according to the differential geometry of the generalized curvature tensor of the almost contact metric manifolds, especially the class of Kenmotsu type and the class of Kenmotsu type as a hypersurface of the Hermitian manifold.

\section{Preliminaries}\label{sec2}
We use the notations $M^{2n+1}$, $X(M)$ and $\nabla$ to denote the smooth manifold $M$ of dimension $2n+1$, the Lie algebra of smooth vector fields of $M$, and the Riemannian connection respectively.
\begin{definition}[\cite{Kir03}]\label{ACR-manifold}
A smooth manifold $M^{2n+1}$ with the quadruple $(\Phi, \xi, \eta, g)$ is called an almost contact metric manifold or briefly ACR-manifold, where $\Phi: X(M)\rightarrow X(M)$, $\xi\in X(M)$, $g$ is the Riemannian metric and $\eta(\cdot)=g(\cdot, \xi)$, are such that
$$\Phi(\xi)=0 ; \quad \eta(\xi)=1 ; \quad \eta\circ\Phi=0 ;\quad \Phi^{2}=-\textrm{id}+\eta\otimes\xi ;$$
$$g(\Phi X, \Phi Y)=g(X, Y)-\eta(X)\eta(Y) ; \quad \forall X, Y\in X(M). $$
\end{definition}
In the present article, we fix the components of the Riemannian metric $g$ of an ACR-manifold $M^{2n+1}$ as follows:
\begin{equation}\label{the metric g}
  g_{00}=1; \quad g_{a0}=g_{ab}=g_{\hat{a}\hat{b}}=0; \quad g_{\hat{a}b}=\delta^{a}_{b}; \quad g_{ij}=g_{ji},
\end{equation}
where $a, b=1, 2, \ldots, n$, $\hat{a}=a+n$ and $i, j=0, 1, \ldots, 2n$. Moreover, the components of the endomorphism $\Phi$ are given by
\begin{equation}\label{the endo. phi}
  \Phi^{0}_{0}=\Phi^{a}_{\hat{b}}=0; \quad \Phi^{a}_{b}=\sqrt{-1}\delta^{a}_{b}; \quad \Phi^{j}_{i}=-\Phi^{\hat{i}}_{\hat{j}},
\end{equation}
where $\hat{\hat{i}}=i$. So, for all $X, Y\in X(M)$, we have
$$ X=X^{i}\varepsilon_{i}; \quad g(X, Y)=g_{ij}X^{i}Y^{j};\quad  \Phi(X)=\Phi^{i}_{j} X^{j}\varepsilon_{i}, $$
where $X^{i}\in C^{\infty}(M)$ and $(p; \varepsilon_{0}=\xi, \varepsilon_{1}, \ldots, \varepsilon_{2n})$ is an A-frame over $M^{2n+1}$ such that $p\in M$, $\varepsilon_{a}=\frac{1}{\sqrt{2}}(id-\sqrt{-1}\Phi)e_{a}$, $\varepsilon_{\hat{a}}=\frac{1}{\sqrt{2}}(id+\sqrt{-1}\Phi)e_{a}$, and $\{\xi, e_{1}, \ldots, e_{n}, \Phi e_{1}, \ldots, \Phi e_{n} \}$ is a basis of $X(M)$. The set of all A-frames as given above is called an associated G-structure space (AG-structure space). For more details, we refer to \cite{Kir03}.
\begin{definition}[\cite{ABA20}]\label{class of Kenmotsu type}
 A class of ACR-manifold such that the following identity:
\begin{equation*}
  \nabla_{X}(\Phi)Y-\nabla_{\Phi X}(\Phi)\Phi Y=-\eta(Y)\Phi X, \quad \forall \; X, Y\in X(M)
\end{equation*}
holds is called a class of Kenmotsu type.
\end{definition}
\begin{lemma}[\cite{ABA20}]\label{smooth functions}
On the AG-structure space, the class of Kenmotsu type satisfies the following relations:
\begin{align*}
A^{ad}_{[bc]}-B^{ad}_{\quad [cb]}-B^{ah}_{\quad [b}\ B^{\quad d}_{|h|c]}&=0; \quad A^{acd}_{b}-B^{a[cd]}_{\quad b}+B^{a[c}_{\quad h} B^{|h|d]}_{\quad b}=0; \quad A^{a}_{[bcd]}=0  \nonumber \\
A_{ad}^{[bc]}+B_{ad}^{\quad [cb]}+B_{ah}^{\quad [b}\ B_{\quad d}^{|h|c]}&=0  ; \quad A_{acd}^{b}+B_{a[cd]}^{\quad b}-B_{a[c}^{\quad h} B_{|h|d]}^{\quad b}=0; \quad A_{a}^{[bcd]}=0;
\end{align*}
where $[\cdot |\cdot| \cdot]$ denotes the anti-symmetric operator of the involved indices except $|\cdot|$ and $c, d, h \in \{1, 2, \ldots, n\}$.
\end{lemma}
We denote by $R$, $r$, $Q$ the Riemann curvature tensor, the Ricci tensor and the Ricci operator of and $ACR$-manifold respectively.
\begin{theorem}[\cite{ABA20}]\label{Riemann components}
The components of $R$ for the class of Kenmotsu type over the AG-structure space are given by
\begin{enumerate}
  \item $R^{a}_{0c0} = -\delta^{a}_{c}; \quad  R^{a}_{\widehat{b}cd} = 2(B^{ab}_{\quad [cd]}-\delta^{a}_{[c}\ \delta^{b}_{d]}); \quad R^{a}_{\widehat{b}c\widehat{d}} = B^{abd}_{\quad c}-B^{ab}_{\quad h}\ B^{hd}_{\quad c};$
  \item $R^{a}_{bcd} = 2A^{a}_{bcd}; \quad R^{a}_{bc\widehat{d}} = A^{ad}_{bc}-B^{ah}_{\quad c}\ B^{\quad d}_{bh} -\delta^{a}_{c}\ \delta^{d}_{b},$
\end{enumerate}
where $R(X, Y)Z=R^{i}_{jkl} X^{k} Y^{l} Z^{j} \varepsilon_{i}$, $k, l=0, 1, \ldots, 2n$ and the remaining components of $R$ are given by the first Bianchi identity or by the conjugate (i.e.\ $R^{i}_{jkl}=\overline{R^{\hat{i}}_{\hat{j}\hat{k}\hat{l}}}$; $\hat{0}=0$) of the above components or are identically zero.
\end{theorem}
\begin{theorem}[\cite{ABA20}]\label{ricci components}
The components of $r$ of the class for Kenmotsu type over the AG-structure space are as follows:
\begin{enumerate}
  \item $r_{00} = -2n; \quad  r_{ab} = -2A^{c}_{abc}+B_{cab}^{\quad c}-B_{ca}^{\quad h}\ B_{hb}^{\quad c};$
  \item $r_{a0} = 0; \quad  r_{\widehat{a}b} = -2(n\delta^{a}_{b}+B^{ca}_{\quad [bc]})+A^{ac}_{cb}-B^{ah}_{\quad b}\ B^{\quad c}_{ch},$
\end{enumerate}
where $r(X, Y)=r_{ij} X^{i} Y^{j}$, $r_{ij}=r_{ji}$ and the remaining components of $r$ are conjugate to the above components.
\end{theorem}
\begin{definition}[\cite{ABA20}]\label{Einstein and phi-invariant Ricci}
An ACR-manifold $(M^{2n+1}, \Phi, \xi, \eta, g)$ with Ricci tensor $r$, 
\begin{enumerate}
  \item is called an Einstein manifold, if $r_{ij}=\lambda g_{ij}$, where $\lambda$ is an Einstein constant.
  \item is called an $\eta$-Einstein manifold, if $r_{ij}=\lambda g_{ij}+\mu\eta_{i}\eta_{j}$, where $\lambda$, $\mu$ are scalars.
  \item is said to have $\Phi$-invariant property, if $r_{a0}=r_{ab}=0$.
\end{enumerate}
\end{definition}
\begin{definition}[\cite{SHK14}]\label{special tensors}
The projective, concircular and generalized curvature tensors of type (4, 0) on the  ACR-manifold $(M^{2n+1}, \Phi, \xi, \eta, g)$ are defined by the following formulas respectively:
\begin{align*}
  P(X, Y, Z, W) & =R(X, Y, Z, W)-\frac{1}{2n}\{ g(X, Z) r(Y, W)-g(X, W) r(Y, Z) \}; \\
  C(X, Y, Z, W) & =R(X, Y, Z, W)-\frac{s}{2n(2n+1)}\{ g(X, Z) g(Y, W)-g(X, W) g(Y, Z) \};\\
  B(X, Y, Z, W) & =a_{0}R(X, Y, Z, W) +a_{1}\{ g(X, Z) r(Y, W)-g(X, W) r(Y, Z) \\
  & \qquad\qquad\qquad\qquad\qquad\quad+r(X, Z)g(Y, W)-r(X, W) g(Y, Z) \} \\
  &{ }\quad + 2a_{2}s\{ g(X, Z) g(Y, W)-g(X, W) g(Y, Z) \};
\end{align*}
for all $X, Y, Z, W\in X(M)$, where $s$ is the scalar curvature, $a_{0}, a_{1}, a_{2}$ are scalars and for any tensor $T$ of type $(3, 1)$, we get
$T(X, Y, Z, W)=g(T(Z, W)Y, X)$, which is a tensor of type $(4, 0)$.
\end{definition}
We can rewrite the above tensors on $AG$-structure space as follows:
\begin{equation}\label{P on AG-structure}
  P_{ijkl}=R_{ijkl}-\frac{1}{2n}\{ g_{ik}\ r_{jl}-g_{il}\ r_{jk} \};
\end{equation}
\begin{equation}\label{C on AG-structure}
  C_{ijkl} =R_{ijkl}-\frac{s}{2n(2n+1)}\{ g_{ik}\ g_{jl}-g_{il}\ g_{jk} \};
\end{equation}
\begin{equation}\label{B on AG-structure}
  B_{ijkl}=a_{0}R_{ijkl}+a_{1}\{ g_{ik}\ r_{jl}-g_{il}\ r_{jk}+r_{ik}\ g_{jl}-r_{il}\ g_{jk} \}+2a_{2}s\{ g_{ik}\ g_{jl}-g_{il}\ g_{jk} \}.
\end{equation}
We note that the generalized curvature tensor $B$ satisfies the first Bianchi identity.

\section{Properties of the Generalized Curvature Tensor}\label{sec3}
In this section, we shall investigate some properties of the generalized curvature tensor on the class of Kenmotsu type.
\begin{theorem}\label{components of B}
On the AG-structure space, the components of the generalized curvature tensor are given by
\begin{enumerate}
  \item $B_{a0b0}=a_{1}\ r_{ab}$;
  \item $B_{\hat{a}0b0}=-(a_{0}+2na_{1}-2a_{2}s)\delta^{a}_{b}+a_{1}\ r_{\hat{a}b}$;
  \item $B_{\hat{a}bcd}=2a_{0}\ A^{a}_{bcd}+a_{1}\{ \delta^{a}_{c}\ r_{bd}-\delta^{a}_{d}\ r_{bc} \}$;
  \item $B_{\hat{a}bc\hat{d}}=a_{0}(A^{ad}_{bc}-B^{ah}_{\quad c}\ B^{\quad d}_{bh})+a_{1}\{\delta^{a}_{c}\ Q_{b}^{d}+\delta^{d}_{b}\ Q_{c}^{a} \}+(2a_{2}s-a_{0}) \delta^{a}_{c}\ \delta^{d}_{b} $;
  \item $B_{\hat{a}\hat{b}cd}=2a_{0}\ B^{ab}_{\quad [cd]}+4a_{1}\ \delta^{[a}_{[c}\ Q_{d]}^{b]}+2(2a_{2}s-a_{0})\ \delta^{[a}_{[c}\ \delta^{b]}_{d]} $;
\end{enumerate}
and the remaining components are identically zero, given by the first Bianchi identity or conjugate to the above components.
\end{theorem}
\begin{proof}
Since $r(X, Y)=g(X, QY)$, then $r_{ij}=g_{ik}Q^{k}_{j}$. Consquently, regarding the Equation~\eqref{the metric g}, we have  $$r_{\hat{a}b}=g_{\hat{a}k}Q^{k}_{b}=g_{\hat{a}0}Q^{0}_{b}+g_{\hat{a}c}Q^{c}_{b}+g_{\hat{a}\hat{c}}Q^{\hat{c}}_{b}=Q^{a}_{b}.$$
Since $B$ is defined on the class of Kenmotsu type, then we may substitute the values of $R_{ijkl}=R^{\hat{i}}_{jkl}$ according to Theorem \ref{Riemann components} and the values of $g_{ij}$ according to Equation~\eqref{the metric g} in Equation~\eqref{B on AG-structure}, obtaining the desired result.
\end{proof}

\begin{theorem}\label{Einstein and B}
The class of Kenmotsu type $(M^{2n+1}, \Phi, \xi, \eta, g)$ has flat generalized curvature tensor if and only if $M$ is $\eta$-Einstein manifold with:
\begin{gather*}
\lambda = \frac{1}{a_{1}} (a_{0}+2na_{1}-2a_{2}s), \quad
A^{a}_{bcd}=0, \quad
\mu = -(2n+\lambda), \\
A^{ad}_{bc} = B^{ah}_{c} B^{d}_{bh} + \frac{a_{1}}{a_{0}} \mu \delta^{a}_{c} \delta^{d}_{b} \quad \textup{ and } \quad
B^{ab}_{[cd]} = \frac{a_{1}}{a_{0}} \mu \delta^{a}_{[c}\delta^{b}_{d]},
\end{gather*}
provided that $a_{0}, a_{1}\neq 0$.
\end{theorem}
\begin{proof}
Suppose that $M^{2n+1}$ has flat generalized curvature tensor with $a_{0}\neq 0$ and $a_{1}\neq 0$, then $B_{ijkl}=0$ and from Theorem \ref{components of B}, we have
$$r_{ab}=0; \quad r_{\hat{a}b}=\frac{1}{a_{1}}(a_{0}+2na_{1}-2a_{2}s)\delta^{a}_{b}; \quad A^{a}_{bcd}=0.$$
Then, according to the Definition \ref{Einstein and phi-invariant Ricci}, we get $\lambda=\frac{1}{a_{1}}(a_{0}+2na_{1}-2a_{2}s)$. Since $M$ is the class of Kenmotsu type, then from the Theorem \ref{ricci components}, we have $r_{00}=-2n=\lambda+\mu$ and this gives us $\mu$. Again, Theorem \ref{components of B}, item 4 gives $A^{ad}_{bc}=B^{ah}_{\quad c}\ B^{\quad d}_{bh}+\frac{a_{1}}{a_{0}}\mu\ \delta^{a}_{c}\ \delta^{d}_{b}$. Moreover, Theorem \ref{components of B}, item 5 gives $B^{ab}_{\quad [cd]}=\frac{a_{1}}{a_{0}}\mu\ \delta^{a}_{[c}\ \delta^{b}_{d]}$. The converse is also true.
\end{proof}
Now, we introduce the notion of generalized $\Phi$-holomorphic sectional $G\Phi HS$-curvature tensor as follows:
\begin{definition}\label{G phi-Holomorphic}
A G$\Phi$HS-curvature tensor $S$ of an ACR-manifold $(M^{2n+1}, \Phi, \xi, \eta, g)$ is a map defined by
$$S(X)=\frac{B(\Phi X, X, X, \Phi X)}{(g(X, X))^{2}}; \quad \forall\ X\in \ker(\eta); \quad X\neq 0.$$
Moreover, $M$ is called of pointwise constant G$\Phi$HS-curvature if $S(X)=\gamma$ and $\gamma$ does not depend on $X$.
\end{definition}
Clearly, a $G\Phi HS$-curvature tensor is a $\Phi$-holomorphic sectional ($\Phi HS$-)curvature tensor if and only if $a_{0}=1$ and $a_{1}=a_{2}=0$. Therefore, we can derive the necessary and sufficient condition for an $ACR$-manifold to have pointwise constant $G\Phi HS$-curvature on $AG$-structure space.
\begin{theorem}\label{G-phi-H-S with ness-suff-cond.}
An ACR-manifold $(M^{2n+1}, \Phi, \xi, \eta, g)$ has pointwise constant G$\Phi$HS-cur\-va\-ture if and only if, on the AG-structure space, the generalized curvature tensor $B$ of $M$  satisfies
$$B^{(a\; d)}_{(bc)}=\frac{\gamma}{2}\widetilde{\delta}^{ad}_{bc},$$
where $\widetilde{\delta}^{ad}_{bc}=\delta^{a}_{b} \delta^{d}_{c}+\delta^{a}_{c} \delta^{d}_{b}$ and $(\cdot \cdot)$ denotes the symmetric operator of the included indices.
\end{theorem}
\begin{proof}
Since the tensor $B$ has the same properties as the Riemannian curvature tensor $R$, then we can follow the same proof was found in \cite{Kir03} or equivalently in \cite{Umn02}.
\end{proof}
\begin{theorem}\label{G-phi-H-S with Kenmotsu type}
\hspace*{-3pt}The class of Kenmotsu type $(M^{2n+1}, \Phi, \xi, \eta, g)$ has pointwise constant G$\Phi$HS\--cur\-va\-ture if and only if, on the AG-structure space, $M$ satisfies the following equality:
$$A^{ad}_{bc}=B^{\quad [ad]}_{bc}-B^{\quad a}_{hb}\ B^{dh}_{\quad c}-\frac{2a_{1}}{a_{0}}\delta^{(a}_{(b}Q^{d)}_{c)}+\frac{\gamma-2a_{2}s+a_{0}}{2a_{0}}\widetilde{\delta}^{ad}_{bc}.$$
\end{theorem}
\begin{proof}
Suppose that $M$ is the class of Kenmotsu type and has pointwise constant G$\Phi$HS-curvature. Using Theorem \ref{G-phi-H-S with ness-suff-cond.} and Theorem \ref{components of B}, item 4, we get
$$A^{(ad)}_{(bc)}=B^{(a|h|}_{\quad (b}\ B^{\quad d)}_{c)h}-\frac{2a_{1}}{a_{0}}\delta^{(a}_{(b}Q^{d)}_{c)}+\frac{\gamma-2a_{2}s+a_{0}}{2a_{0}}\widetilde{\delta}^{ad}_{bc}.$$
The above equation can be rewritten as follows:
$$A^{(ad)}_{(bc)}=-B^{\quad (a}_{h(b}\ B^{d)h}_{\quad c)}-\frac{2a_{1}}{a_{0}}\delta^{(a}_{(b}Q^{d)}_{c)}+\frac{\gamma-2a_{2}s+a_{0}}{2a_{0}}\widetilde{\delta}^{ad}_{bc}.$$
Since $A^{ad}_{bc}=A^{[ad]}_{[bc]}+A^{[ad]}_{(bc)}+A^{(ad)}_{[bc]}+A^{(ad)}_{(bc)}$, then taking into account Lemma \ref{smooth functions} and the above result, we conclude the proof.
\end{proof}
Recently, Y{\i}ld{\i}z and De \cite{YID12} introduced the notions of $\Phi$-projectively semisymmetric and $\Phi$-Weyl semisymmetric. Regarding these ideas, we can introduce the following definition:
\begin{definition}\label{phi-G-semisymmetric}
An ACR-manifold $(M^{2n+1}, \Phi, \xi, \eta, g)$ is called $\Phi$-generalized semisymmetric if $B(Z, W)\cdot\Phi=0$, for all $Z, W\in X(M)$, or equivalently
$$B(X, \Phi Y, Z, W)+B(\Phi X, Y, Z, W)=0; \quad  \forall\ X, Y, Z, W\in X(M).$$
\end{definition}
\begin{lemma}\label{phi-G-semisymmetric on AG-space}
On AG-structure space, the ACR-manifold $(M^{2n+1}, \Phi, \xi, \eta, g)$ is $\Phi$-generalized semi ($\Phi$GS-)symmetric if and only if
$$B_{a0b0}=B_{\hat{a}0b0}=B_{a0bc}=B_{\widehat{a}0bc}=B_{a0\widehat{b}c}=B_{abcd}=B_{\hat{a}\hat{b}cd}=0.$$
\end{lemma}
\begin{proof}
According to the Definition \ref{phi-G-semisymmetric}, we have that $M$ is $\Phi$-generalized semi-symmetric if and only if
$$B(X, \Phi Y, Z, W)+B(\Phi X, Y, Z, W)=0; \quad  \forall\ X, Y, Z, W\in X(M).$$
On the $AG$-structure space, the above identity is equivalent to the following:
$$B_{iqkl}\ \Phi^{q}_{j}+B_{tjkl}\ \Phi^{t}_{i}=0; \quad  q, t=0, 1, \ldots, 2n.$$
If we take $$(i, j, k, l)=(a, 0, b, 0), (\hat{a}, 0, b, 0), (a, 0, b, c), (\widehat{a}, 0, b, c), (a, 0, \widehat{b}, c), (a, b, c, d), (\hat{a}, \hat{b}, c, d),$$ and using Equation~\eqref{the endo. phi}, we obtain the result.
\end{proof}
It is not hard to conclude the following:
\begin{corollary}\label{phi-G-semisymmetric with flat B}
The ACR-manifold $(M^{2n+1}, \Phi, \xi, \eta, g)$ of flat generalized curvature tensor is usually $\Phi$GS-symmetric.
\end{corollary}
\begin{corollary}\label{phi-GS-symmetric with K-type and flat B}
The class of Kenmotsu type $(M^{2n+1}, \Phi, \xi, \eta, g)$ has flat generalized curvature tensor if and only if $M$ is $\Phi$GS-symmetric with $A^{a}_{bcd}=0$ and $A^{ad}_{bc}=B^{ah}_{\quad c}\ B^{\quad d}_{bh}+\frac{a_{1}}{a_{0}}\mu\ \delta^{a}_{c}\delta^{d}_{b}$, where $\mu=-\frac{1}{a_{1}}(a_{0}+4na_{1}-2a_{2}s)$, provided that $a_{0}, a_{1}\neq 0$.
\end{corollary}
\begin{proof}
Suppose that $M$ is the class of Kenmotsu type and it has flat generalized curvature tensor, then from the Corollary \ref{phi-G-semisymmetric with flat B}, $M$ is $\Phi$GS-symmetric and from the Theorem \ref{components of B}, we get the other conditions.

Conversely, If $M$ is $\Phi GS$-symmetric with the above conditions then according to Lemma~\ref{phi-G-semisymmetric on AG-space} and Theorem \ref{components of B}, $M$ has flat generalized curvature tensor.
\end{proof}
\begin{theorem}\label{phi-G-semisymmetric with K-type}
The class of Kenmotsu type $(M^{2n+1}, \Phi, \xi, \eta, g)$ is $\Phi$GS-symmetric if and only if $M$ is $\eta$-Einstein manifold with $\lambda=\frac{1}{a_{1}}(a_{0}+2na_{1}-2a_{2}s)$, $\mu=-(2n+\lambda)$ and $B^{ab}_{\quad [cd]}=\frac{a_{1}}{a_{0}}\mu\ \delta^{a}_{[c}\delta^{b}_{d]}$, provided that $a_{0}, a_{1}\neq 0$.
\end{theorem}
\begin{proof}
Suppose that $M$ is $\Phi$-generalized semi-symmetric class of Kenmotsu type, then from Lemma \ref{phi-G-semisymmetric on AG-space} and Theorem \ref{components of B}, we have
$$r_{ab}=0; \quad r_{\hat{a}b}=\frac{1}{a_{1}}(a_{0}+2na_{1}-2a_{2}s)\delta^{a}_{b}; \quad B^{ab}_{\quad [cd]}=-\frac{1}{a_{0}}(a_{0}+4na_{1}-2a_{2}s) \delta^{a}_{[c}\delta^{b}_{d]}.$$
Regarding the Definition \ref{Einstein and phi-invariant Ricci} and Theorem \ref{ricci components}, we obtain the values of $\lambda$ and $\mu$.

The converse is verified directly.
\end{proof}
\begin{corollary}\label{phi-GS-symmetric with G-phi-HS-curvature}
The class of Kenmotsu type $(M^{2n+1}, \Phi, \xi, \eta, g)$ is $\Phi$GS-symmetric and has G$\Phi$HS-curvature if and only if, $M$ is $\eta$-Einstein manifold with $\lambda=\frac{1}{a_{1}}(a_{0}+2na_{1}-2a_{2}s)$, $\mu=-(2n+\lambda)$, $B^{ab}_{\quad [cd]}=\frac{a_{1}}{a_{0}}\mu\delta^{a}_{[c}\delta^{b}_{d]}$, and $A^{ad}_{bc}=\frac{\gamma}{2a_{0}}\widetilde{\delta}^{ad}_{bc}-B^{\quad a}_{hb}\ B^{dh}_{\quad c}+\frac{a_{1}}{a_{0}}\mu\delta^{a}_{b}\delta^{d}_{c}$, provided that $a_{0}, a_{1}\neq 0$.
\end{corollary}
\begin{proof}
Suppose that $M$ is the class of Kenmotsu type, then the necessary and sufficient conditions of the present corollary are satisfied from Theorems \ref{G-phi-H-S with Kenmotsu type} and \ref{phi-G-semisymmetric with K-type}.
\end{proof}

\section{Generalized Curvature Tensor Related with Another Tensors}\label{sec4}
In this section, we introduce a generalization of the notion of $ACR$-manifold of constant curvature that is used by Abood and Al-Hussaini \cite{ABH19}. We shall present this idea in the following definition:
\begin{definition}\label{constant generalized curvature def.}
An ACR-manifold $(M^{2n+1}, \Phi, \xi, \eta, g)$ is said to have constant generalized curvature $\kappa$ if the following identity holds:
$$B(X, Y, Z, W)=\kappa\{g(X, Z)g(Y, W)-g(X, W)g(Y, Z)\}; \quad \forall\ X, Y, Z, W\in X(M).$$
\end{definition}
On the $AG$-structure space, Definition \ref{constant generalized curvature def.} equivalent to the identity below.
\begin{equation}\label{constant generalized curvature on AG-space}
  B_{ijkl}=\kappa\{g_{ik}\ g_{jl}-g_{il}\ g_{jk}\}.
\end{equation}
Directly, regarding the Definition \ref{constant generalized curvature def.}, Definition \ref{special tensors} and the definition of the conharmonic curvature tensor (see \cite{DES19}), we have the following result:
\begin{theorem}\label{constant generalized curvature with conharmonic}
Suppose that $M^{2n+1}$ is an ACR-manifold of constant generalized curvature $\kappa=2a_{2}s$. Then $M$ has flat conharmonic curvature tensor if and only if, $a_{0}=1$ and $a_{1}=-\frac{1}{2n-1}$.
\end{theorem}
\begin{theorem}\label{cons. gen. curv. components}
An ACR-manifold $(M^{2n+1}, \Phi, \xi, \eta, g)$ has constant generalized curvature $\kappa$ if and only if, on the AG-structure space, $B$ has the following components:
\begin{enumerate}
  \item $B_{\hat{a}0b0}=\kappa\ \delta^{a}_{b}$;
  \item $B_{\hat{a}bc\hat{d}}=\kappa\ \delta^{a}_{c}\delta^{d}_{b}$;
  \item $B_{\hat{a}\hat{b}cd}=2\kappa\ \delta^{a}_{[c}\delta^{b}_{d]}$;
\end{enumerate}
and the remaining components are identically zero, obtained from the above components by the first Bianchi identity or by taking the conjugate operation.
\end{theorem}
\begin{proof}
The result follows from Equation~\eqref{constant generalized curvature on AG-space} by taking
$$(i, j, k, l)=(\hat{a}, 0, b, 0), (\hat{a}, b, c, \hat{d}), (\hat{a}, \hat{b}, c, d);$$
and using the Equation~\eqref{the metric g}.
\end{proof}
\begin{theorem}\label{phi-GS-symmetric with G.C.C.}
The ACR-manifold $(M^{2n+1}, \Phi, \xi, \eta, g)$ is $\Phi$GS-symmetric if and only if, $M$ has constant generalized curvature $\kappa=0$.
\end{theorem}
\begin{proof}
The claim of this theorem is obtained from Lemma \ref{phi-G-semisymmetric on AG-space} and Theorem \ref{cons. gen. curv. components}.
\end{proof}
\begin{theorem}\label{G-phi-H-S with C.G.C.}
If an ACR-manifold $(M^{2n+1}, \Phi, \xi, \eta, g)$ has constant generalized curvature $\kappa$, then $M$ has pointwise constant G$\Phi$HS-curvature equal to $\gamma=\kappa$.
\end{theorem}
\begin{proof}
The result follows from Theorems \ref{G-phi-H-S with ness-suff-cond.} and \ref{cons. gen. curv. components}.
\end{proof}
\begin{theorem}\label{C.G.C. with Kenmotsu type}
The class of Kenmotsu type $(M^{2n+1}, \Phi, \xi, \eta, g)$ has constant generalized curvature $\kappa$ if and only if, $M$ is an $\eta$-Einstein manifold with:
\begin{gather*}
\lambda=\frac{1}{a_{1}}(a_{0}+2na_{1}-2a_{2}s+\kappa), \quad
A^{a}_{bcd}=0, \quad
\mu=-(2n+\lambda), \\
A^{ad}_{bc} = B^{ah}_{c} B^{d}_{bh} + \frac{a_{1}}{a_{0}} \mu \delta^{a}_{c} \delta^{d}_{b} \quad \textup{and} \quad
B^{ab}_{[cd]} = \frac{a_{1}}{a_{0}} \mu \delta^{a}_{[c} \delta^{b}_{d]},
\end{gather*}
provided that $a_{0}, a_{1}\neq 0$.
\end{theorem}
\begin{proof}
The assertion of this theorem is obtained by combining the results of the Theorems~\ref{components of B} and \ref{cons. gen. curv. components}.
\end{proof}
Now, we find the geometric properties of $ACR$-manifold if the generalized curvature tensor, the concircular tensor and the projective tensor are related.

Suppose that $(M^{2n+1}, \Phi, \xi, \eta, g)$ is an $ACR$-manifold satisfies the following condition:
\begin{equation}\label{dependent tensors}
  B(X, Y, Z, W)=\frac{a_{0}}{3}\{P(X, Y, Z, W)-P(Y, X, Z, W)+C(X, Y, Z, W)\}.
\end{equation}
Regarding the Equations~\eqref{P on AG-structure}, \eqref{C on AG-structure} and \eqref{B on AG-structure}, we can write the Equation~\eqref{dependent tensors} on the $AG$-structure space as follows:
\begin{equation}\label{dependent tensors on AG-structure}
  (a_{1}+\frac{a_{0}}{6n})\{ g_{ik}\ r_{jl}-g_{il}\ r_{jk}+r_{ik}\ g_{jl}-r_{il}\ g_{jk} \}+(2a_{2}+\frac{a_{0}}{6n(2n+1)})s\{ g_{ik}\ g_{jl}-g_{il}\ g_{jk} \}=0.
\end{equation}
The contracting of the Equation~\eqref{dependent tensors on AG-structure}, that is, multiplying it by $g^{ik}$ (the components of $g^{-1}$ on $AG$-structure space), we can deduce that
\begin{equation}\label{contract dependent tensors}
   r_{jl}=-\frac{(\alpha+2n\beta)s}{(2n-1)\alpha} g_{jl},
\end{equation}
where $\alpha=a_{1}+\frac{a_{0}}{6n}$ and $\beta=2a_{2}+\frac{a_{0}}{6n(2n+1)}$. Moreover, the contracting of Equation~\eqref{contract dependent tensors} gives $a_{0}+4na_{1}+4n(2n+1)a_{2}=0$. Then we can state the following theorem:
\begin{theorem}\label{dependent tensors with Einstein}
Any ACR-manifold $(M^{2n+1}, \Phi, \xi, \eta, g)$ which satisfies the identity \eqref{dependent tensors} is an Einstein manifold with $a_{0}+4na_{1}+4n(2n+1)a_{2}=0$, provided that $\alpha\neq 0$. Moreover, if $M$ is the class of Kenmotsu type then $s=\frac{2n(2n-1)\alpha}{\alpha+2n\beta}$, provided that $\alpha+2n\beta\neq 0$.
\end{theorem}
\begin{proof}
The first part of this theorem is obvious from the above discussion. Now, if $M$ is the class of Kenmotsu type then from the Theorem \ref{ricci components}, we have $r_{00}=-2n$. Then the result is established from the Equations~\eqref{the metric g} and \eqref{contract dependent tensors}.
\end{proof}

\section{The Hypersurfaces of the Hermitian Manifold}\label{sec5}
Suppose that $(M^{2n-1}, \Phi, \xi, \eta, g)$ is an $ACR$-manifold, then there exists an almost complex structure $J$ on $M\times\mathbb{R}$ defined by $J(X, f\frac{d}{dt})=(\Phi X-f\xi, \eta(X)\frac{d}{dt})$, where $X\in X(M)$, $t\in\mathbb{R}$ and $f$ is a smooth function on $\mathbb{R}$. The Riemannian metric $h$ on $M\times\mathbb{R}$ is defined by
$$h((X, f_{1}\frac{d}{dt}), (Y, f_{2}\frac{d}{dt}))=g(X, Y)+f_{1}\ f_{2};\quad \forall\ X, Y\in X(M);\quad f_{1}, f_{2}\in C^{\infty}(\mathbb{R}).$$
The structure on $M\times\mathbb{R}$ is Hermitian if and only if the structure on $M$ is normal (see \cite{Pit07}). Since the class of Kenmotsu type is normal because it is the class $C_{3}\oplus C_{4}\oplus C_{5}$, where $C_{5}$ is taken here to be $\alpha$-Kenmotsu manifold with $\alpha=1$ (see \cite{CHG90} for more detail about the classes $C_{3}$ and $C_{4}$). Then the product manifold of the class of Kenmotsu type and the real line is Hermitian (i.e.\ $W_{3}\oplus W_{4}$, see \cite{GRH80}).

Now, we discuss the opposite problem, that is, if $(N^{2n}, J, h)$ is an Hermitian manifold, then can we find a hypersuface of $N$ which is the class of Kenmotsu type? We rely on the citation \cite{BAK15} for the background.

Suppose that $a, b, c=1, 2, \ldots, n-1$ and $\sigma_{ij}=\sigma_{ji}$; $i, j=1, 2, \ldots, 2n-1$ are the components of the second quadratic form as mentioned in \cite{BAK15}.
\begin{theorem}[\cite{BAK15}]\label{hypers. thm}
An ACR-manifold which a hypersuface of an almost Hermitian manifold has the following first family of the Cartan structure equations:
\begin{align*}
  d\omega^{a} & =\omega^{a}_{b}\wedge\omega^{b}+B^{ab}_{c}\ \omega^{c}\wedge\omega_{b}+B^{abc}\ \omega_{b}\wedge\omega_{c}+(\sqrt{2}B^{an}_{b}+\sqrt{-1}\sigma^{a}_{b})\omega^{b}\wedge\omega \\
   & +(\sqrt{-1}\sigma^{ab}-\sqrt{2}\widetilde{B}^{nab}-\frac{1}{\sqrt{2}}B^{ab}_{n}-\frac{1}{\sqrt{2}}\widetilde{B}^{abn})\omega_{b}\wedge\omega;\\
 d\omega_{a} & =-\omega_{a}^{b}\wedge\omega_{b}+B_{ab}^{c}\ \omega_{c}\wedge\omega^{b}+B_{abc}\ \omega^{b}\wedge\omega^{c}+(\sqrt{2}B_{an}^{b}-\sqrt{-1}\sigma_{a}^{b})\omega_{b}\wedge\omega \\
   & -(\sqrt{-1}\sigma_{ab}+\sqrt{2}\widetilde{B}_{nab}+\frac{1}{\sqrt{2}}B_{ab}^{n}+\frac{1}{\sqrt{2}}\widetilde{B}_{abn})\omega^{b}\wedge\omega;\\
d\omega &=\sqrt{2}B_{nab}\ \omega^{a}\wedge\omega^{b}+\sqrt{2}B^{nab}\ \omega_{a}\wedge\omega_{b}+(\sqrt{2}B^{na}_{b}-\sqrt{2}B_{nb}^{a}-2\sqrt{-1}\sigma^{a}_{b})\omega^{b}\wedge\omega_{a}\\
&+(\widetilde{B}_{nbn}+B^{n}_{nb}+\sqrt{-1}\sigma_{nb})\omega\wedge\omega^{b}+(\widetilde{B}^{nbn}+B_{n}^{nb}-\sqrt{-1}\sigma_{n}^{b})\omega\wedge\omega_{b}.
\end{align*}
\end{theorem}
From Banaru \cite{Ban15}, we see that the Hermitian manifold $N$ satisfies $B^{\alpha\beta\gamma}=B_{\alpha\beta\gamma}=0$, where $\alpha, \beta, \gamma=1, 2, \ldots, n$, and thus the Theorem \ref{hypers. thm} reduce to the following form:
\begin{theorem}\label{hypers. Hermitian}
An ACR-manifold which a hypersuface of the Hermitian manifold has the following first family of the Cartan structure equations:
\begin{align*}
  d\omega^{a} & =\omega^{a}_{b}\wedge\omega^{b}+B^{ab}_{c}\ \omega^{c}\wedge\omega_{b}+(\sqrt{2}B^{an}_{b}+\sqrt{-1}\sigma^{a}_{b})\omega^{b}\wedge\omega
  +(\sqrt{-1}\sigma^{ab}-\frac{1}{\sqrt{2}}B^{ab}_{n})\omega_{b}\wedge\omega;\\
 d\omega_{a} & =-\omega_{a}^{b}\wedge\omega_{b}+B_{ab}^{c}\ \omega_{c}\wedge\omega^{b}+(\sqrt{2}B_{an}^{b}-\sqrt{-1}\sigma_{a}^{b})\omega_{b}\wedge\omega
 -(\sqrt{-1}\sigma_{ab}+\frac{1}{\sqrt{2}}B_{ab}^{n})\omega^{b}\wedge\omega;\\
d\omega &=(\sqrt{2}B^{na}_{b}-\sqrt{2}B_{nb}^{a}-2\sqrt{-1}\sigma^{a}_{b})\omega^{b}\wedge\omega_{a}+(B^{n}_{nb}
+\sqrt{-1}\sigma_{nb})\omega\wedge\omega^{b}\\
&+(B_{n}^{nb}-\sqrt{-1}\sigma_{n}^{b})\omega\wedge\omega_{b}.
\end{align*}
\end{theorem}
Regarding Abood and Abass \cite{ABA20}, we note that the class of Kenmotsu type satisfies the following theorem:
\begin{theorem}[\cite{ABA20}]\label{stq}
The class of Kenmotsu type $M^{2n-1}$ has the following first group of Cartan structure equations:
\begin{align*}
  d\omega^{a}&=\omega^{a}_{b}\wedge\omega^{b}+ B^{ab}_{\quad c}\ \omega^{c}\wedge\omega_{b}-\omega^{a}\wedge\omega;\\
  d\omega_{a}&=-\omega^{b}_{a}\wedge\omega_{b}+ B^{\quad c}_{ab}\ \omega_{c}\wedge\omega^{b}-\omega_{a}\wedge\omega;\\
  d\omega&=0,
\end{align*}
where $B^{ab}_{\quad c}$ and $B^{\quad c}_{ab}$ are the components of the first Kirichenko's tensor as explained in \emph{\cite{KID06}}.
\end{theorem}
Now, if the class of Kenmotsu type $M^{2n-1}$ is a hypersurface of the Hermitian manifold $N^{2n}$, then the cartan structure equations that mentioned in the Theorems \ref{hypers. Hermitian} and \ref{stq} must be equal. Then we get
\begin{align}\label{equal hypers.}
  B^{ab}_{c} & =B^{ab}_{\quad c}; \quad  \sqrt{2}B^{an}_{b}+\sqrt{-1}\sigma^{a}_{b}=-\delta^{a}_{b}; \quad \sqrt{-1}\sigma^{ab}-\frac{1}{\sqrt{2}}B^{ab}_{n}=0; \nonumber \\
  B_{ab}^{c} & =B^{\quad c}_{ab}; \quad \sqrt{2}B_{an}^{b}-\sqrt{-1}\sigma_{a}^{b}=-\delta^{b}_{a}; \quad \sqrt{-1}\sigma_{ab}+\frac{1}{\sqrt{2}}B_{ab}^{n}=0;\\
  \sqrt{2}B^{na}_{b}&-\sqrt{2}B_{nb}^{a}-2\sqrt{-1}\sigma^{a}_{b}=0; \quad B^{n}_{nb}
+\sqrt{-1}\sigma_{nb}=0; \quad B_{n}^{nb}-\sqrt{-1}\sigma_{n}^{b}=0. \nonumber
\end{align}
Since $\sigma_{[\alpha\beta]}=0$ and $B_{[\alpha\beta]}^{\gamma}=B_{\alpha\beta}^{\gamma}$, then Equation~\eqref{equal hypers.} gives the following relations:
\begin{equation}\label{second form}
  \sigma_{ab}=0; \quad \sigma_{nb}=0; \quad \sigma^{a}_{b}=\sqrt{-1}(\sqrt{2}B^{an}_{b}+\delta^{a}_{b}).
\end{equation}
Then from the above discussion, we can establishing the theorem below.
\begin{theorem}\label{Kenmotsu with Hermitian}
If the Hermitian manifold has the class of Kenmotsu type as a hypersurface, then the second quadratic form has components agree with the Equation~\eqref{second form}.
\end{theorem}
On the other hand, we can establish a relation between the components of the Riemannian curvature tensors of the almost Hermitian manifold and its hypersurfaces.

Suppose that $\mathcal{R}^{i}_{jkl}$ are the components of the Riemannian curvature tensor of the almost Hermitian manifold, $N^{2n}$ and $\mathcal{\widetilde{R}}^{i}_{jkl}$ are the components of the Riemannian curvature tensor of its hypersurface $M^{2n-1}$. Then from the second group of cartan structure equations, we have
\begin{align*}
  d\omega^{i}_{j} & =\omega^{i}_{k}\wedge\omega^{k}_{j}+\frac{1}{2}\mathcal{R}^{i}_{jkl}\ \omega^{k}\wedge\omega^{l}; \\
  d\theta^{i}_{j} & =\theta^{i}_{k}\wedge\theta^{k}_{j}+\frac{1}{2}\mathcal{\widetilde{R}}^{i}_{jkl}\ \theta^{k}\wedge\theta^{l};
\end{align*}
where $\omega^{i}_{j}$ and $\theta^{i}_{j}$ are the Riemannian connection forms of $N$ and $M$ respectively. Whereas, $\omega^{k}$ and $\theta^{k}$ are the dual $A$-frames on $AG$-structure spaces of $N$ and $M$ respectively. Moreover, from \cite{BAK15}, we have
$$\theta^{i}=C^{i}_{j}\ \omega^{j}; \quad \omega^{i}=\widetilde{C}^{i}_{j}\ \theta^{j}; \quad \theta^{i}_{j}=C^{i}_{k}\ \omega^{k}_{r}\ \widetilde{C}^{r}_{j}; \quad \omega^{i}_{j}=\widetilde{C}^{i}_{k}\ \theta^{k}_{r}\ C^{r}_{j};$$
where $C=(C^{i}_{j})$ and $C^{-1}=(\widetilde{C}^{i}_{j})$ were defined in \cite{BAK15}. Then the substitution of the above relations in the second group of cartan structure equations, we conclude the following theorem:
\begin{theorem}\label{last thm}
If $\mathcal{R}^{i}_{jkl}$ and $\mathcal{\widetilde{R}}^{q}_{rst}$ are the components of the Riemannian curvature tensor of the almost Hermitian manifold $(N^{2n}, J, g)$ and its hypersurface $(M^{2n-1}, \Phi, \xi, \eta, g)$ respectively, then they are related as follow:
$$\mathcal{R}^{i}_{jkl}=\widetilde{C}^{i}_{q}\ \mathcal{\widetilde{R}}^{q}_{rst}\ C^{r}_{j}\ C^{s}_{k}\ C^{t}_{l}.$$
\end{theorem}


\EditInfo{June 28, 2021}{December 31, 2021}{Haizhong Li}

\end{paper}